\documentclass[12pt, leqno]{amsart}

\usepackage{amsfonts}
\usepackage{amssymb}
\usepackage{amsmath}
\usepackage{amscd}
\usepackage{float}
\usepackage{graphicx}
\usepackage{amsfonts}

\usepackage{color}
\usepackage{cancel}%\sout{};  \cancel{}
\usepackage[cp1250]{inputenc}
\usepackage[T1]{fontenc}

\usepackage{multirow}
\usepackage{multicol}
\usepackage{hyperref}

\usepackage{color}
\usepackage{xypic}

\usepackage[cmtip,arrow]{xy}

\usepackage[normalem]{ulem}
\newtheorem{theorem}{Theorem}

\newtheorem{lemma}{Lemma}
\newtheorem{proposition}{Proposition}
\newtheorem{definition}{Definition}
\newtheorem{remark}{Remark}
\newtheorem*{exa}{\bf Example}

\def\A{\mathsf{A}}
\def\AA{\mathcal{A}}
\def\B{\mathcal{B}}

\def\D{\mathsf{D}}

\def\b{\mathsf{b}}

\def\dd{\mathsf{d}}
\def\DD{\mathsf{D}}

\begin{document}

\title{Catalan  Triangles  and  Tied Arc Diagrams}

\author{F. Aicardi}
 \address{Sistiana 56, Trieste IT}
 \email{francescaicardi22@gmail.com}

\begin{abstract} The  Catalan triangle, as  well  as  a   Fuss-Catalan triangle,   enter a  problem  of  counting particular    tied  arc  diagrams.   This  setting  allows us to  prove  some  combinatorial properties of  these  triangles.
\end{abstract}



\date{}
\maketitle
\section{Introduction}
%While Catalan  numbers $C_n =\frac{1}{n+1}\binom{2n} {n}$ enter a huge set of  different  %combinatorial problems,  the  recurrence  for  them  here introduced  is perhaps  new.
The  Catalan  numbers $C_n :=\frac{1}{n+1}\binom{2n} {n}$ can be   decomposed  into  $n$  integers  defined by a  recurrence: they form  the so-called Catalan triangle (see \cite{Rog,Sha}).  Other  decompositions  exist for the  Fuss-Catalan  numbers $\A_n(p,q):=  \frac{q}{pn+q}\binom{pn+q}{n}$ (see
\cite{Fin}).

Since  the seminal  paper  \cite{aijuJKTR1}, several  tied  knot  algebras  of  different  kinds  were  introduced  and  studied.  In particular, for some families   of  such algebras,  a  diagrammatical interpretation   of  the  generators   leads to the definition of  tied  arc  diagrams.

The  tied  arc  diagrams  can  be  obtained  recursively. It turns out that their number, in the  case  of  the so  called tb-diagrams,  is calculated by means of the  Catalan  triangle. In  the case of  the so called ta-diagrams the  recurrence yields  a Fuss-Catalan triangle, whose  rows  sum  to  the Fuss-Catalan numbers $\A_n(4,1)$.
  We use  the  tb-diagrams  to  prove  a  combinatorial  result (Theorem \ref{t2}), that is  used in \cite{AAJ} to calculate  the  dimension
 of  a tied-tangle  monoid.

 The  combinatorial  results of Theorems \ref{t3}, proved by  means of  the  ta-diagrams, are  used in    \cite{AJP}
 to get the  dimensions  of  tied Temperley-Lieb algebras.
\section{Results}

\subsection{The  Catalan Triangle}

\begin{definition}\cite{Rog,Sha}\label{def1}\rm The  integers $T(n,k)$  are  defined for  $n\ge0$,  $0\le k \le n$  by the initial conditions,
\begin{equation}\label{tr1}
T(0,0)=1, \quad
  T(n,0)=1,\quad  T(n,n)=0, \quad  \quad n>0 ;
  \end{equation}
  and  the  recurrence
 \begin{equation} \label{tr2}  T(n,k)= T(n,k-1)+T(n-1, k), \quad   \quad  n>1, \quad 0<k < n .
 \end{equation}
 \end{definition}
Observe  that  the  above equations  imply
  \begin{equation} \label{tr3}  T(n,k)=  \sum_{j=0}^k T(n-1,j) \quad \quad  n>1, \quad 0<k < n.
 \end{equation}
Here  the  Catalan  triangle $T(n,k)$, for  $n\le 8$,  $0\le k\le n$
{\tiny
 $$\begin{array}{ccccccccccccccccccccccccc}
          n& &   &  &   &   &   &   &   &   &   &     &   &   &   &    &   &   &   &   &     &  &     \\
          0& &   &  &   &   &   &   &   &   &   &     & 1  &   &   &    &   &   &   &   &     &  &    \\
          1& &   &  &   &   &   &   &   &   &   &  1  &    & 0 &   &    &   &   &   &   &     &  &    \\
          2& &   &  &   &   &   &   &   &   & 1 &     & 1  &   & 0 &    &   &   &   &   &     &  &    \\
          3& &   &  &   &   &   &   &   & 1 &   & 2   &    & 2 &   & 0  &   &   &   &   &     &  &    \\
          4& &   &  &   &   &   &   & 1 &   & 3 &     & 5  &   & 5 &    & 0 &   &   &   &     &  &     \\
          5& &   &  &   &   &   & 1 &   & 4 &   & 9   &    & 14&    & 14 &   & 0 &   &   &     &  &   \\
          6& &   &  &   &   & 1 &   & 5 &   & 14 &    & 28 &   & 42 &   & 42 &   &  0 &   &     &  &    \\
          7& &   &  &   & 1 &   & 6 &   & 20 &   & 48 &    & 90&    & 132&   & 132&   & 0  &     &  &    \\
          8& &   &  & 1 &   & 7 &   &27 &    & 75&    &165 &   & 297&   & 429 &    &429 &    & 0   &  &    \\
         \end{array}$$
}
 \begin{proposition}\label{t1}\cite{Rog} The  Catalan  triangle  satisfies, for every $n>0$ :
 $$ T(n,n-1)= C_{n-1},  \quad \sum_{k=0}^n  T(n,k)=C_{n}. $$
 \end{proposition}

\begin{theorem}\label{t2} The  Catalan  triangle  satisfies, for every $n>0$ :
 $$ \sum_{k=0}^n  T(n,k) 2^{n-1-k}=\binom{2n-1}{n}. $$
 \end{theorem}

%%%%%%%%%%%%%%%%%%%%%%%%%%%%%%%%%%%%%%%%%%%%%%%%%%%%%%%%%%%
\subsection{A  Fuss-Catalan triangle   } $\quad$

 The triangle here  defined by  recurrence was  already obtained in  \cite{Fin}  by  another  procedure.

\begin{definition}\label{def2}\rm The  integers $F(n,k)$  are  defined for  $n\ge 0$,  $0\le k \le n$  by the initial conditions
\begin{equation}\label{E1}
F(0,0)=1; \quad
  F(n,n)=0, \quad \text{ for} \quad  n>0 ;
  \end{equation}
  and  the  recurrence for $n>0$  and  $0\le  k<n$:
 \begin{equation} \label{E2}   F(n,k)=\sum_{j=0}^{k}  \binom{k-j+2}{2} F(n-1, j).
 \end{equation}
 \end{definition}

Here  the  triangle $F(n,k)$  for  $n\le 8$,  $0\le k\le n$
{\tiny
 $$\begin{array}{ccccccccccccccccccc}
          n&   &   &   &   &   &   &   &   &     &   &       &   &    &   &   &   &   &   \\
          0&   &   &   &   &   &   &   &   &     & 1  &   &   &    &   &   &   &   &   \\
          1&   &   &   &   &   &   &   &   &  1  &    & 0 &   &    &   &   &   &   &   \\
          2&   &   &   &   &   &   &   & 1 &     & 3 &   & 0 &    &   &   &   &   &   \\
          3&   &   &   &   &   &   & 1 &   & 6   &    & 15&   & 0  &   &   &   &   &   \\
          4&   &   &   &   &   & 1 &   & 9 &     & 39 &   & 91&    & 0 &   &   &   &   \\
          5&   &   &   &   & 1 &   & 12&   & 72  &    & 272&    & 612 &   & 0 &   &   &   \\
          6&   &   &   & 1 &   & 15 &   & 114 &    & 570 &   & 1995 &   & 4389 &   &  0 &   &   \\
          7&   &   & 1 &   & 18 &   & 165 &   & 1012 &    & 4554&    & 15180&   & 32890&   & 0  &   \\
          8&   & 1 &   & 21 &   &225&    & 1625&    &8775&   & 36855&   &118755 &    &254475 &    & 0 \\
           \end{array}$$
}
\begin{theorem}  \label{t3}
The  integers $F(n,k)$    satisfy,  for  every  $n>0$
\begin{equation}\label{f3}  \sum_{k=0}^{n} F(n,k)= \A_n(4,1);  \end{equation}

\begin{equation}\label{f4}  \sum_{k=0}^{n} F(n,k)\binom{n-k+3}{3}= \A_n(4,4);  \end{equation}
  and
\begin{equation} \label{f5}  F(n,n-1)=\A_{n-1}(4,3)  \end{equation}
where    $\A_n(4,1 )$, $\A_n(4,4 )$  and $\A_n(4,3) $  are  respectively   the  Fuss-Catalan  numbers  $\frac{1}{4n+1} {{4n+1}\choose{n}}$,   $\frac{4}{4n+4} {{4n+4}\choose{n}}$ and  $\frac{3}{4n+3} {{4n+3}\choose{n}}$.
\end{theorem}

\begin{remark}\rm  Statements \ref{f3}  and  \ref{f5} can be obtained  also  from   \cite{Fin}. However,  we give  here  a proof in terms of  ta-diagrams for  completeness.
\end{remark}

\section{  Proof  of  Proposition 1 in terms of arc diagrams}
 Proposition \ref{t1}  says  that  the  numbers $T(n,k)$  form a  particular partition  of $C_n$.
We  will  consider
 one  of  the definitions  of  $C_n$,  namely:
 {\it $C_n$ is  the  number of  semicircle  diagrams  with $n$  semicircles}.
More precisely:
\begin{definition} \rm An {\it $n$-arcdiagram }  is  a planar  diagram  of $n$  non intersecting semicircles or  {\it  arcs}  with  end points on a straight line. All  arcs    lie  on one of the  closed  half--planes  defined  by the line.
\end{definition}

\begin{lemma}\label{l1} The  number  of $n$-arcdiagrams is the  Catalan  number  $C_n$\end{lemma}
\begin{proof} Recall  that
 the  Catalan  number  $C_n$   is, among others,  the  number of expressions containing $n$ pairs of parentheses which are correctly matched.  E.g., for  $n=3$, $C_n=5$:
$$   ( \quad )( \quad )( \quad ), \quad ((\quad))(\quad), \quad (\quad)((\quad )), \quad ((\quad)(\quad)), \quad (((\quad))).   $$
To  obtain  the   corresponding  five arcdiagrams,  it is  sufficient  to  substitute  every  pair  of open-closed parentheses  with  a semicircle in this way:
$$   (\quad)  \quad  \rightarrow \quad  (\  \includegraphics[scale=1.2,trim=0.4cm  0 0.4cm 0 ]{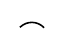} \ ) \quad \rightarrow   \includegraphics[scale=1.5,trim=0 0.4cm 0 0]{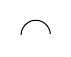} $$
This defines   evidently a  bijection,  since  from an $n$-arcdiagram  we get viceversa a  unique set of $n$ pairs of parentheses.
\end{proof}
Observe  that  an  $n$-arcdiagram  divides the  upper half plane  in $n+1$  regions,  of  which  only one,  said $U$,  is unbounded.

 \begin{definition} \rm  A {\sl block} of  an  arcdiagram, is a semidisc  defined  by  an  arc at  the  boundary  of  $U$. A block   may contain   other  arcs.
 \end{definition}

 Let us denote $D_n^k$  the  number  of the   $n$-arcdiagrams having    $k$  blocks.

 \begin{lemma}\label{l2} $D_n^k=T(n,n-k)$. \end{lemma}
 \begin{proof}
 Evidently, $D_n^n=1$, i.e., there is  one  diagram  with $n$ blocks, each one  with a  sole   semicircle,  and  $D_n^0=0$,  since a  diagram  with no  blocks  has  no  arcs;  so, in particular,  $D_0^0=1$. The initial conditions (\ref{tr1})  are fulfilled.
 We  will  verify  the  recurrence  (\ref{tr3})  that  reads
 \begin{equation}\label{Dr}   D_n^k= \sum_{j=k-1}^{n-1} D_{n-1}^j. \end{equation}
    Consider  the  set  of  arcdiagrams with $n-1$  arcs. It is  partitioned  into  $n-1$ parts corresponding  to  the  number  $j$ of  blocks, so  that  $C_{n-1}= \sum_{j=1}^{n-1}D_{n-1}^j$.
Suppose  that we  want  to  get,  starting  from this  set,   the
set of   arcdiagrams with $n$  arcs that have  exactly  $k$ blocks, by  inserting  in a  suitable  way  the  $n$-th  arc. Let  $\dd$  be  one of  such  diagrams. If  the  number $j$ of  blocks of  $\dd$  is  less
than  $k-1$,  the  new  diagram cannot  have  $k$ blocks  by  adding a  new  arc. Therefore  we  start  with  diagrams with $j\ge k-1$.    We   put the  left endpoint of  the  $n-th$  arc at  left of the diagram $\dd$.  If  $\dd$ has  $k-1$ blocks,  then  the  right point of the new arc  will  be at  left  of  the  first  block of $\dd$.  If  $\dd$  has  $j>k-1$ arcs, the end  point  will  be  after  the  $j-k+1$th block,  thill  $j=n-1$. To  finish,   we  have  to  prove  that in this  way  we  get  all  diagrams with  $n$  arcs  and  $k$  blocks,  and that  all  diagrams  so  obtained  are  all  different.  Firstly,  for  every  $j$  we  consider a  set  of diagrams with  $j$  blocks and  $n-1$ arcs,  that  are  all  different,  and hence the  procedure we  use  produces  different diagrams.  Two diagrams  obtained  starting  from two sets  with  different values of $j$  cannot  coincide  since  the  first  blocks of  them  contain  different  quantities of  blocks of  the  original  diagrams. Suppose now  that a  diagram  $\dd'$  with  $k$  blocks  and  $n$  arcs  is  not  reached  by  the  above procedure.   The  diagram   $\dd'$    contains one first block. Removing  the  first arc, we get  a sequence (possibly  empty) of $h$ blocks $\b_1,\dots,\b_h$    $0\le h\le n-k$, plus  $k-1$  blocks $\b_{h+1},\dots, \b_{h+k}$. The  sequence   of blocks $\b_1,\dots \b_{h+k-1}$ form a  diagram  with  $n-1$  arcs  and  with  a  number of  blocks  at  least  $k-1$. Then  it  is impossible that it has been missed by  the  procedure.
\end{proof}
\begin{proof}[Proof  of Proposition 1]
By  Lemma  \ref{l2}, we have  to  verify  that  $D_{n}^1=C_{n-1}$. Indeed, this is  what  says Eq. (\ref{Dr}) with  $k=1$,  by  using Lemma \ref{l1}.
As for  the  second  equation,  we  have by Lemma \ref{l1}, $  \sum_{k=1}^n  D_n^k= C_n$, that  we can rewrite as
 $$  \sum_{k=0}^{n-1}  D_{n}^{n-k}= C_n.$$
 Then Proposition \ref{t1}  follows from Lemma  \ref{l2}, since $D_n^0=0$.
\end{proof}

%\begin{remark}  $C_n$ is the  number of  non-isomorphic plane  trees with  $n+1$  vertices.
% The set of $n$-sc--diagrams  is in  bijection  with  the  set   $T_{n}$ of     ordered trees with $n+1$ vertices. (An %ordered tree is a rooted tree in which the children of each vertex are given a fixed left-to-right order).   The  %bijection  is  shown  in the  following picture, where  $n=6$.
%\end{remark}

%\centerline{\includegraphics[scale=0.8]{tree1.pdf}}

\section{Tb-diagrams  and proof of  Theorem \ref{t2}}

\begin{figure}  [H]
 \centering
 \includegraphics{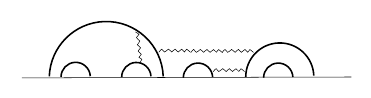}
 \caption{A tied arcdiagram}\label{F1}

\end{figure}

 \rm  A  {\it tied} arcdiagram  is an    arcdiagram  with  {\it ties} that  may  connect each other  two  arcs,  avoiding   selfintersections  and intersections  with the  arcs. See  Figure \ref{F1}.

%Particular  tied  arc  diagrams  are  concerned  with  some  tied  versions  of  the  Temperly-Lieb  algebra.

We  firstly  consider  particular tied  arcdiagrams named  tb-diagrams.

 \begin{definition}\label{toparc}\rm We  call  {\sl top arc} of  an arcdiagram  $\dd$ a  semicircle  bounding a  block  of  $\dd$. A  diagram  has  at least one top  arc. The  top  arcs are  naturally  ordered from left to right.\end{definition}

 \begin{definition}\label{tb}\rm
 A  {\sl  tb-diagram}  is a  tied arcdiagram   whose  ties may  exist only in  the  unbounded  region of  the  complement to the diagram and can   connect  only successive   top arcs.
\end{definition}

\begin{lemma}  \label{l3} The  number  of  tb-diagrams  with $n$  arcs  and  $k$ blocks is  $2^{k-1}D_n^k$.
\end{lemma}
\begin{proof} It  is  an  immediate  consequence  of Definition \ref{tb},  since there  are  $k-1$ pairs of  successive   top arcs  admitting a  tie  in between.
\end{proof}

 We will define  a  bijection  between the   tb-diagrams  with  $n$  arcs,  and  the  $n$-combinations  of $2n-1$  objects.  Then Theorem \ref{t2} will  follow  from Lemma \ref{l2}  and Lemma \ref{l3}.

%We  show a     procedure  to  associate to any  tb-diagram a unique  $n$-combination of  $\{1,2,\dots,2n-1\}$, and  a  %second  procedure  to  associate to  any  $n$-combination a unique  tb-diagram.

  Let  $\dd$  be  a  tb-diagram  with  $n$ arcs. Label   by  $0, 1, \dots, {2n-1}$  the  endpoints  of  its  arcs.

We  associate  to $\dd$  one  $n$-combination of $\{1,2,\dots,2n-1\}$  as follows.

{\it Procedure 1.}
\begin{itemize}
\item If  $\dd$ has no  ties, then  take  the $n$ labels of the  right endpoints   of  the  $n$  arcs.
\item If  $\dd$ has ties,  then  take  the   labels of the   the right endpoints  of  the  arcs of the  first block,  and   the  left  ones of the arcs inside  every  block connected  to  the preceding block  by a  tie.
\end{itemize}

\begin{figure}
 [H]
 \centering
 \includegraphics{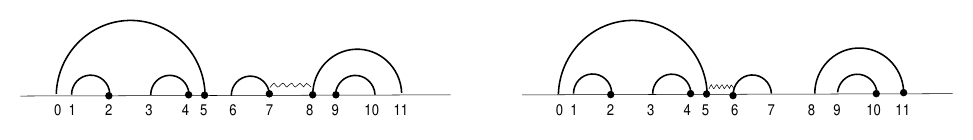}
 \caption{The   6-combinations  are (2,4,5,7,8,9)  and (2,4,5,6,10,11)}\label{F2}
\end{figure}

\begin{exa}\rm In Figure \ref{F2} see two  tb-diagrams  with 6  arcs, and  the  corresponding  combinations  of  the integers from 1 to  11.
\end{exa}

\begin{lemma} Every  $n$-combination of $\{1,2,\dots,2n-1\}$  is  obtained by    Procedure 1  starting  from  one  and only one  tb-diagram.
\end{lemma}

\begin{proof} Consider  firstly  a  tb-diagram  without  ties  and  a  function  $f$ defined  on the  points $i$, taking  value -1  on  the left  endpoint and +1  on the  right  endpoint of each  arc. Define the function $F$ as
$$     F(p):=\sum_{j=0}^{p}f(j)$$
Let  $\{2k_i-1 \}_{i=1}^m$  be  the zeroes of  $F$.  Now  we  observe  the trivial  facts:
 \begin{enumerate}
 \item  $f$  has  at  least one  zero  for  $k_m=n$;
   \item the  points $ {2k_i-1}$  are the  right endpoints of the  top  arcs of  $\dd$,  so $\dd$ has  $m$  blocks;
 \item   the point  $p_0$  and   and the points   $p_{2k_i}$ are the  left  endpoints  of  the  top  arcs;
 \end{enumerate}
 Consider  the arcdiagrams   inside  any  block of  $\dd$,  obtained by  removing the top  arc.   The function $f$  takes  value $ -1$  at  the  first  points  $p_{2k_i+1}$  and $+1$  at  the  last  points $ 2k_{i+1}$.
    Then,  by  relabeling  the  indices   $(2k_i+1, \dots  2k_{i+1})$ by  $(0,\dots ,2n-1)$,  we  can  repeat all  the  preceding  observations,  and  so  on  for  every  subdiagram inside  the  blocks  of  the  arcdiagrams  just  considered.

 It is therefore  evident  that $f$  defines  uniquely the   arcdiagram,  and  that  the integers  $i$    such  that  $f(i)=+1$  define one  $n$-combination  of  $\{1,2\dots, 2n-1\}$.

 It is  also  clear  that  the $n$-combinations   obtained  this  way are  not all:  for  instance, all  such  combinations  contain   the  integer $2n-1$,  since  for  every arcdiagram $f({2n-1})=+1$.

 We observe  now  that the  function $-f$  defines  the  same arcdiagram  as $f$,  simply  exchanging left endpoints with  right  endpoints. Moreover,  every  function $f'$ obtained  from $f$ by  reversing the   values   inside  one or  more  blocks, still  defines  the  same arcdiagram as $f$:  the  value +1  will be  assigned  to  the  left endpoints  inside  the  blocks  where     $f'=-f$.
\begin{figure}
 [H]
 \centering
 \includegraphics{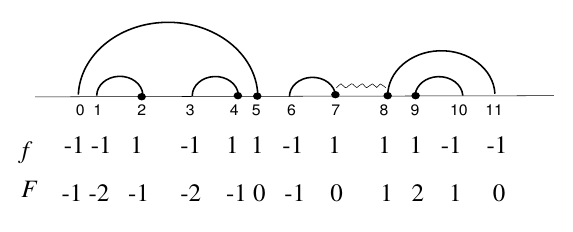}
 \caption{The functions $f$  and  $F$  of  the  left  tb-diagram of  Figure \ref{F2}}\label{F3}
\end{figure}

 Thus Procedure 1 associates   different  $n$-combinations to all  $2^{k-1}$  tb-diagrams obtained  by  adding  ties  to  an arcdiagram $\dd$ with  $k$ blocks.  By  the  preceding observations,    the  functions $f$  associated  to  these  tb-diagrams  differ  each  other only  by  a  reversion of   sign  inside  some   blocks of  $\dd$.

Now we  assign     uniquely one  tb-diagram   to  a chosen  $n$-combination.

{\it Procedure 2.}

Let us  start  with a $n$-combination  $C$ of $\{1,2,\dots 2n-1\}$.
  Consider  $2n$ points on a  line labeled  by
 $\{0,1,\dots, {2n-1}\}$  and   a  function $f$  defined on these  points:
   $$  f(i)=+1 \quad  \text{if} \  i\in C, \quad  f(i)=-1 \quad \text{otherwise} \  $$
    Observe  that  $f(0)=-1$.
 Define  $F(p)=\sum_{i=0}^p f(i)$. Since  the  quantity  of  zeroes  and  ones  is  the  same,  surely  $F(p)=0$ at least  at  $p=2n-1$.  Suppose  there  are  no  other zeroes.  It  means  that $F$ is  negative till  $i={2n-2}$,  and  that $f({2n-1})=+1$.  I.e., $2n-1 \in C$. So, $\dd$  has  only one  block,  and  $(0,2n-1)$ are  the  endpoints of    the unique  top arc.  Remove  this top  arc  and define  $F^1(k)=\sum_{i=1}^k f(i)$. Note that $f(1)=-1$,  otherwise $1$  should be  a  zero  of  $F$.   The  zeroes   of  $F^1$  defines  the right  endpoints  of  the  top arcs inside  the  block of  $\dd$. Observe  that  $F^1(2n-2)=0$,  and  if  $k$ is  a  zero  of $F^1$, then $f({k+1})=-1$,
 otherwise  $ k+1$ should  be  a  zero  of  $F$.  We  proceed  this  way  inside the  blocks defined  by  the  top  arcs, determining  all new top  arcs, and  so  on  inside  all corresponding  blocks,  concluding   that  at  the left endpoint  $l$ of every  arc,  $f(l)=-1$,  since   every  arc  is  the  top  arc  of one  block  of  some  subdiagram  of  $\dd$.  The  diagram  $\dd$  is thus  completely defined. Observe  that in this  case $\dd$ and has no  ties.

 Suppose  now  that  $F(i)=0$  for  $i=2m-1$,  $m<n$.  The  points  $(0,{2m-1})$  define  the top arc  of  the  first  block  of  $\dd$,  and  inside  it  we  proceed  as  in the  preceding  case,   putting  all  arcs, that  result to  have  left  endpoints  where  $f=-1$.  Observe  that this  follows  only from the fact    that  $f(0)=-1$  by  hypothesis.

 We  proceed now   to  define  the  second  block of  $\dd$.  If $f(2m)=-1$,  and  $F(2r-1)=0$ for  $m<r\le n$,  then $f(2r-1)=+1$,  and we  proceed as  previously  inside  the  second  block.  But   if $f(2m)=+1$,  and  $F(2r-1)=0$ for  $m<r\le n$,  then $f(2r-1)=-1$,  so  that  on the  endpoints of the second  top arc $(2m,2r-1)$,  $f $  takes  value +1 on  the  left,  and  -1 on the  right.  Then  we  proceed  as  in  the  preceding  case  but  exchanging +1  with  -1, so that the  second  block  is uniquely  defined,  and  the  left endpoints of all  arcs  inside  it   have  $f=+1$.  In this  case  we put a  tie  between  the  first  and  the  second  block  of  $\dd$.

 We  proceed  the  same  way  to  define  the  successive  blocks of  $\dd$,  tied or  not  with  the preceding one.
The  tb-diagram $\dd$ is  uniquely  defined  by  the $n$-combination.

 It is  clear  that, having obtained  by Procedure 2 a  tb-diagram  $\dd$  from a chosen   $n$-combination, the same  combination  is  obtained  by applying to  $\dd$  Procedure 1.  So  we  have proved  that by Procedure  1  all  $n$-combinations  are attained.

\end{proof}

\section{Ta-diagrams    }

In this  section  we introduce  another  class  of  tied  arcdiagrams,  called ta-diagrams.

In fact,   ta-diagrams  are  equivalence  classes  of  tied arcdiagrams,  the  equivalence  relation  being given below.

We order the  arcs  of a  $n$-arcdiagram   according to the  order of  their right  endpoints on  the  line  from  left  to  right,  see  Figure \ref{F4}.

In a ta-diagram with   arcs  $a_1,\dots,a_n$,  the  ties define  a  partition  of the  set  of  arcs in this  way: {\sl two arcs connected by a  tie  belong to  the  same part of  the partition.} We shall denote this partition  by ta-partition,  and  for  short  we  write only the  indices $1,\dots,n$ of  the  arcs.

\begin{definition} \rm Two ta-diagrams are  equivalent if they coincide forgetting the ties, and  the  ta-partition  defined by  the  ties  of  one   diagram coincides with that of  the other.
\end{definition}
\begin{exa} \rm See  Figure \ref{F4}.
\begin{figure}  [H]
 \centering
 \includegraphics{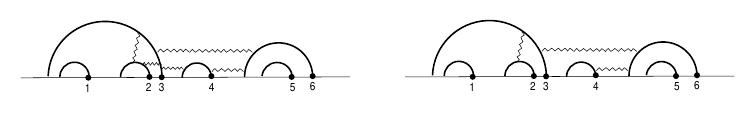}
 \caption{Two equivalent ta-diagrams  with ta-partition $\{\{1 \},\{2,3,4,6\},\{5 \}\}$}\label{F4}
 \end{figure}
 \end{exa}

We  say  that  $m>2$  ties  form a cycle, if  there is a  sequence of $m$ arcs  $a_{i_1},\dots,a_{i_m}$ such  that  there is a  tie between $a_{i_k}$ and  $a_{i_{k+1}}$, for $k<m$, and there is a tie   between $a_{i_m}$ and $a_{i_1}$.

\begin{definition} \rm A  ta-diagram  is  said {\it irreducible}  if  between  two  arcs there is at most one tie,  and  the  ties    do not form cycles.  Otherwise, the ta-diagram is  said  {\it  reducible}.
 \end{definition}

 \begin{proposition} Every ta-diagram   is equivalent to  an irreducible  one.
\end{proposition}

\begin{proof} If  a ta-diagram is reducible,  then  it  becomes irreducible by   canceling   all  ties  between two  arcs  but one, and by canceling  just one  tie, for  every  cycle  of $m$ arcs connected  by  $m$ ties. By these canceling operations the ta-partition defined by the  ties is  preserved.
\end{proof}

\begin{exa} \rm The  left ta-diagram of  Figure \ref{F4} is  reducible, the  right diagram is irreducible.
\end{exa}

%\begin{remark}\rm \label{trees} An irreducible $n$--ta-diagram with $p$ parts,  has $n-p$ ties. Since in each part the ties do not form  cycles,  for each part
%the  arcs  are the  vertices of a tree,  whose edges are the  ties.
%\end{remark}

Observe that an arc $a$  divides the half--plane in two parts  that  we  call the {\it interior}  and  the {\it exterior} of $a$.

\begin{definition}\rm  An  arc $a$ is {\it inside} an  arc $b$  if  it lies in the interior of $b$, otherwise it is {\it outside} $b$.
\end{definition}

\begin{definition}\label{standard} \rm  A ta-diagram is  {\it standard} if it is irreducible and  in each part  consisting of  $m$ arcs  with ordered  indices  $i_1<i_2< \dots < i_m$,
 every arc with  index  $i_k$ is connected by a tie at most to  one arc of index $i_r > i_k$,  and  at most to two  arcs of indices $i_p<i_q< i_k$,   with   $a_{i_p}$ outside  $a_{i_k}$,  and  $a_{i_q}$ inside   $a_{i_k}$.
\end{definition}

\begin{exa}\rm The irreducible  ta-diagram in  Figure  \ref{F4} is not standard, since  the arc $a_6$ is connected by ties to  $a_3$ and  $a_4$, both outside $a_6$. The corresponding standard diagram  is in Figure \ref{F5}, at left.

\begin{figure}  [H]
 \centering
 \includegraphics{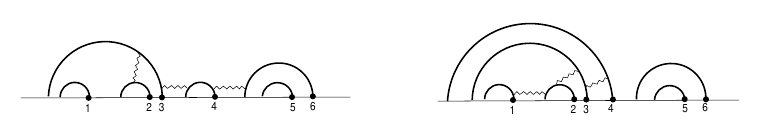}
 \caption{Standard ta-diagrams  }\label{F5}
 \end{figure}

\end{exa}
%\begin{definition}     The {\it  level } of  an  arc is zero  if its internal is void.    The level of an arc is $k$ %if its internal contains  a $(k-1)$-level arc.
%\end{definition}

\begin{proposition}\label{pro1} Every ta-diagram is  equivalent to a  standard ta-diagram.
\end{proposition}
%To each part of the partition  with $m$ arcs,  there  correspond  $m-1$ ties.

\begin{proof} Consider an irreducible  ta-diagram.    Suppose  that   the  diagram is not standard.  Then,  let the arc $a_h$ be the arc with the minimum index for     which  the conditions of being standard are not fulfilled. There are three possibilities:

1) $a_h$   is connected  by ties  to two or more arcs inside it;

2) $a_h$  is connected  by ties to  two or more arcs with higher indices;

3)  $a_h$  is connected  by ties to  two or more arcs with lower  indices outside it.

We  see  how  these   situations  can be corrected to obtain  a standard  diagram,  without affecting the ta-partition.

Observe that in the case (1), there are $m\ge 2$  arcs $a_{i_1}, \dots a_{i_m}$ inside $a_h$  tied with $a_h$.   The arcs  are ordered by  their  increasing indices  $i_1< \dots < i_m$. (in Figure \ref{F6} $m=2$). We  cancel  $m-1$  ties  between  $a_{i_j}$ and $a_h$  for $j<m$    and  we put a tie  from  $a_{i_j}$ and $a_{i_{j+1}}$  for  $i=1, \dots,m-1$.
In the  case (2), we order as before  the   $m\ge 2$ arcs  connected to $a_h$, we cancel  the  ties  but  that from the first  arc to  $a_h$  and  we put a tie  from  the other successive  arcs,  so that it remains only   the tie  from   $a_h$ to the arc with minimum index.

In the case (3) we do the same as in  case  (2).

Now,  observe that the part of the ta-partition  containing $a_h$  remains   unchanged  by  the  above  operations and the ta-diagram obtained  is     irreducible.
It remains to prove  that   we can  always put  the new  ties  without  crossing other arcs  or other  ties.

Consider  the  case (1).   The arcs  $a_{i_j}$  and  $a_{i_{j+1}}$  were initially connected by a  tie to the arc $a_h$. See  Figure  \ref{F6}. This means  that  both  these  arcs were not  inside other arcs between them  and $a_h$.  Then  a  tie  between  $a_{i_j}$  and  $a_{i_{j+1}}$   does not cross any  arc.  This  tie  could  cross  another  tie  between        an arc $a_k$   with  $i_j<k<i_{j+1}$  and  another  arc  with  index  $l$  with  $l<i_j$  or $i_{j+1}<l<h$,  belonging to  another part of the ta-partition. (In  Figure  \ref{F6}    hypothetical ties  of such kind are shown as dotted lines.)  But such a  tie does not  exist  since  it should cross the  tie  connecting  $a_h$  either to $a_{i_j}$  or to $a_{i_{j+1}}$,  see Figure \ref{F6}.

\begin{figure}  [H]
 \centering
 \includegraphics{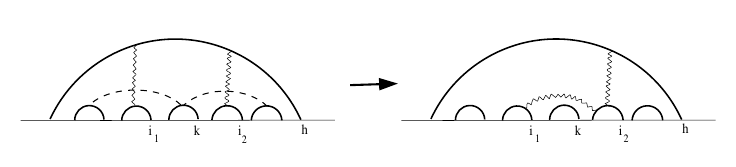}
 \caption{Creating a new tie  between the arcs $a_{i_1}$  and  $a_{i_{2}}$ }\label{F6}
 \end{figure}

The  cases  (2)  and (3)  are  similar.

\end{proof}

\begin{definition} \rm
A {\it block} of a ta-diagram,  is a  subdiagram, having either only one top arc, or a number $m>1$ of top arcs such that the  first   and the  the $m$-th   belong to the same part of the partition.
\end{definition}

\begin{figure}  [H]
 \centering
 \includegraphics{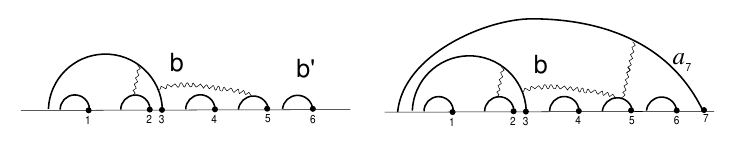}
 \caption{Left. A standard   ta-diagram  $\dd$ with 6 arcs and   two  blocks, $\b$ and $\b'$. Right. A tie  connects the arc  $a_7$    to  the  block $\b$.    }\label{F7}
 \end{figure}

\begin{exa} \rm In  Figure  \ref{F5}, the ta-diagram  at  left  has  a  sole block, the ta-diagram  at  right  has  two blocks. In Figure\ref{F7}, the  ta-diagram  at  left  has  two  blocks.
\end{exa}

In the  next  section  we  will  build the standard  ta-diagrams with $n$  arcs from  those with  $(n-1)$ arcs.  Taking  a  ta-diagram $\dd$ with $n-1$ arcs,  we will   add a  new    arc  $a_n$ and a possible tie from $a_n$ to  some  arc of  $\dd$. Observe that this  arc  is   the top  arc   with maximum index  of a bloc of $\dd$.   We  will  say  for  short  that  the  tie connects  $a_n$  to  that  block.  See  Figure \ref{F7}.

\section{Proof of  Theorem \ref{t3}}

We denote by $\AA_n$  the set of classes of   ta-diagrams with $n$ arcs  and by $A_n$ the  number of standard ta-diagrams. By  Proposition \ref{pro1}, $A_n=|\AA_n|$.

We denote by $\B^j_n$ the  set of classes of  ta-diagrams having $n$  arcs and $j$ blocks,  and $B^j_n$ its cardinality.  Then
\begin{equation}\label{An}  A_n=  \sum_{j=1}^n B_n^j.
\end{equation}

\begin{remark}\rm  Observe that    $\AA_0$  contains  the diagram  without  arcs,  and  $\B_0^0$ the  diagram  with no  arcs and no  blocks,  whereas  $\B_n^0$ is empty, for $n>1$. So
 \begin{equation}\label{iniz} A_0=1,  \quad  B_0^0=1,  \quad   B_n^0=0  \quad \text{for} \quad n>0  \end{equation}
\end{remark}

The proof of Theorem  \ref{t3} is  subdivided  in two  parts.   In  the  first  part  we  prove  that
 \begin{equation}\label{part1} A_n=  \A_n(4,1),  \quad    B^1_n=\A_{n-1}(4,3).  \end{equation}
  In  the  second  part  we  prove  that
  \begin{equation}\label{part2}    F(n,k)=B_n^{n-k} .  \end{equation}

  \subsection{Part 1}  To  prove Equations (\ref{part1}),
we firstly observe  that $A_n$  can be  written  through the  $A_k$  and the $B^1_{n-k}$,  for  $0\le k<n$.

\begin{lemma}\label{l5} For  every $n>0$, the  numbers  $A_k$ and $B^1_{n-k}$, for  $k=0,\dots,n-1$,    satisfy  the equation:
\begin{equation}\label{AnBn}
A_n= \sum_{k=0}^{n-1}   B^1_{n-k} A_k.
\end{equation}
\end{lemma}

\begin{proof}  Consider any standard ta-diagram $\dd$ in $\AA_n$, and decompose it in two ta-diagrams $\DD'$  and $\DD''$:  $\D''$ consists of the last block  of $\dd$,    $\DD'$  is  the  remaining  ta-diagram, which  is void  when $\dd$ has only one  block, i.e. when $\DD''=\dd$. Consider now    the ta-diagrams  for which  $\DD'$ and  $\DD''$  contain  respectively $k$ and $n-k$ arcs.  Their  number  is  equal to  the number of ta-diagrams   with   $k$  arcs,  $A_k$, multiplied by the number of  ta-diagrams with one  block  and $n-k$ arcs, i.e.  $B^1_{n-k}$. Since for  every standard  ta-diagram $\dd$ in $\A_n$  this  decomposition is  unique,  we  get
$$  A_n=  A_0 B_{n}^1 + A_1 B^1_{n-1} +   A_{n-2} B^1_2 +\dots + A_{n-1} B^1_1 $$
i.e.,  eq. (\ref{AnBn}).
\end{proof}

In order  to  use  Equation (\ref{AnBn}) as  a  recurrence  to calculate $A_n$,  we  need  the  following  lemma.

\begin{lemma}  The  number  $B_n^1$ of  standard  ta-diagrams  with $n$  arcs and only one  block
 satisfies
\begin{equation}\label{Bn1}
B_n^1=  \sum_{k=0}^{n-1} A_k \sum_{j=0}^{n-1-k} A_j A_{n-1-k-j}.
\end{equation}
\end{lemma}
\begin{proof} Observe  that   formula (\ref{Bn1})  says that  $B^1_n$   equals    the  sum of all  products of  three ordered  factors $A_iA_jA_k$ such  that  $i+j+k=n-1$. So,  we  prove   that, for any ordered triple of non negative integers $(i,j,k)$     such that $i+j+k=n-1$,   we  can  associate to every triple of standard ta-diagrams in $\AA_i \times \AA_j \times \AA_k$,   one and only  one  standard  diagram  in  $\B^1_n$, and  that each  ta-diagram  in $\B^1_n$ can be  uniquely  decomposed in one triple  of  ta-diagrams in $\AA_i \times \AA_j \times \AA_k$  for some ordered triple $(i,j,k)$ such that $i+j+k=n-1$.

Let us take  three  standard  ta-diagrams,  $\DD_1\in \AA_i$, $\DD_2\in \AA_j$  and $\DD_3\in \AA_k$, with $i+j+k=n-1$.  The  diagram in $\AA_0$ is a diagram  without arcs.   We define a  standard ta-diagram in $B_n^1$   the  following  way: we  put $\DD_2$  at  right of $\DD_1$,  and  $\DD_3$  at  right of $\DD_2$. We get  a  diagram in $\A_{n-1}$  which is  standard.   We  add  one $n$-th  arc, $a_n$, with  right endpoint  at  right of $\DD_3$, and  left endpoint   between  $\DD_1$ and $\DD_2$. Then we put one  tie  from $a_n$  to    the first block of $\DD_1$, if  $\DD_1$ is  not void, and  one  tie  from  $a_n$  to   the  last block of $\DD_2$, if $\DD_2$ is  not void.
Such operations guarantee  that  the new  diagram  is  standard, is uniquely defined starting by $\DD_i,\DD_j,\DD_k$ and has  one  block by  construction.  See  Figure \ref{F8}.

\begin{figure}  [H]
 \centering
 \includegraphics{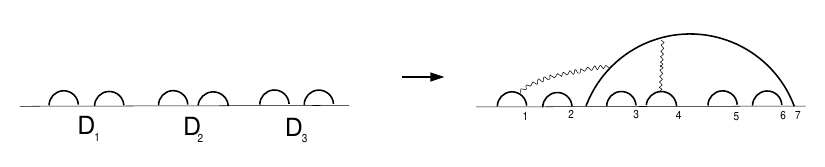}
 \caption{  A standard   ta-diagram      with one block  built from   three  standard  diagrams  $\DD_1,\DD_2,\DD_3\in \AA_2$.  Ties  connect the arc  $a_7$    to  the first  block of $\DD_1$  and   to  the last block  of  $\DD_2$.    }\label{F8}
 \end{figure}

Now,  consider a  standard  diagram  $\dd\in \B_n^1$.  Take the last top arc of  the unique block of  $\dd$, i.e. $a_n$. If the top arc is  unique,  then  $i=0$.  Otherwise, $a_n$ is  connected  by a  tie  to another top  arc. By canceling  this tie,  we  get a  ta-diagram  with  $i<n$ arcs  at  left  of  $a_n$.  We call  this  diagram $\DD_1$.   Since  $\dd$ is  standard,  there is  at most  one  tie from  $a_n$  to  an  arc  $a_m$ inside $a_n$.  If  there is  no such a tie, then $j=0$.  Otherwise,  consider the  ta-diagram  $\DD_2$
consisting of  all  arcs with indices  from $i+1$  to  $m$, i.e.  $j=m-i$.  Observe that  all  these  arcs are  between the left  endpoint of  $a_n$  and  the right endpoint of $a_m$,  since $a_m$  is  connected by a  tie  to  $a_n$.  The  arcs  from $a_{m}$  to  $a_{n-1}$ form  the  ta-diagram  $\DD_3\in \AA_k$,  $k=n-1-m$.   If  $m=n-1$,  then $k=0$.  See  Figure \ref{F9}.
Observe  that  the  diagrams  $\DD_1$, $\DD_2$, $\DD_3$  are  in  this  way  uniquely  defined.

\begin{figure}  [H]
 \centering
 \includegraphics{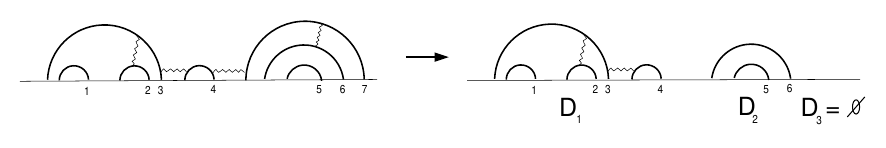}
 \caption{  The standard   ta-diagram  at left      is  decomposed in  $\DD_1\in \AA_4$, $\DD_2\in \AA_2$, $\DD_3\in \AA_0$.    }\label{F9}
 \end{figure}

\end{proof}

Now, we recall  in the  next two lemmas \cite{Riordan}   two known properties of the  Fuss-Catalan  numbers.

\begin{lemma}\label{L1} $\A_n(p,q)=\A_n(p,q-1)+\A_{n-1}(p,p+q-1)$.
\end{lemma}

\begin{lemma} \label{L2} $\A_n(p,r+s)=\sum_{i=0}^n \A_i(p,r)\A_{n-i}(p,s)$.
\end{lemma}

  By  Lemma \ref{L1} for $p=4$, $q=1$  we  have
\begin{equation}\label{formula5}
\A_n(4,1)=\A_n(4,0)+\A_{n-1}(4,4)=\A_{n-1}(4,4),
\end{equation}
since $\A_n(p,0)=0$ for every $p$.

By lemma \ref{L2}  we obtain, for $p=4$, $r=1$ and $s=3$
$$\A_{n-1}(4,4)=\sum_{i=0}^{n-1} \A_i(4,1)\A_{n-1-i}(4,3).$$
Then,  by  Eq. (\ref{formula5})
\begin{equation}\label{formula6}
 \A_{n}(4,1)=\sum_{i=0}^{n-1} \A_i(4,1)\A_{n-1-i}(4,3).
\end{equation}

Using again Lemma \ref{L2}, for $p=4$, $r=1$ and $s=2$ we  get
$$ \A_m(4,3)= \sum_{k=0}^m \A_k(4,1)\A_{m-k}(4,2); $$
and  again for  for $p=4$, $r=1$ and $s=1$,
$$ \A_{m-k}(4,2)=  \sum_{j=1}^{m-k} \A_j(4,1) \A_{m-k-j}(4,1)  . $$
Therefore,
\begin{equation}\label{formula7}
 \A_m(4,3)=\sum_{k=0}^{m} \A_k(4,1) \sum_{j=1}^{m-k} \A_j(4,1) \A_{m-k-j}(4,1).
\end{equation}

Observe  now  that  Equations (\ref{AnBn})  and (\ref{Bn1}),  with  initial  condition (\ref{iniz}),  define  by  recursion the  values $B^1_n$  and  $A_n$,  for  every  integer $n>0$.

Notice, moreover, that  $\A_0(4,1)=1$,  and  $\A_0(4,3)=1$.
Comparing   now  Equations  (\ref{AnBn})  and (\ref{Bn1})  respectively  with   Equations (\ref{formula6}) and (\ref{formula7}),  we  get Equations (\ref{part1}).  So, part 1 is proved.

%At the end  of this section,  we will prove the  following theorem.
%  \begin{theorem} \label{T1} The  number  of standard  ta-diagrams  with  $n$  arcs    is  the Fuss-Catalan  number $C_n(4,1)$, i.e.:
%\begin{equation}    \label{F2}A_n=\frac{1}{4n+1} {{4n+1}\choose{n}}\end{equation}
%\end{theorem}

%%%%%%%%%%%%%%%%%%%%%%%%%%%%%%%%%%%%%%%%%%%%%%%%%%%%%%%%%

\subsection{Part 2.} To prove Equation (\ref{part2}),  we  prove  that  the  integers $B_n^{n-k}$ satisfy  the  same  recurrence of  $F(n,k)$.

I.e.,  for  a given $n$, we  write the integers $ B_n^k$, for $0<k\le n$ in term of the integers $B_{n-1}^j$, for $k-1\le j\le n-1$.

Firstly, we prove this   lemma
\begin{lemma}\label{L9} We  have
\begin{equation}  A_n=  \sum_{k=1}^{n-1}  {{k+3}\choose{3}} B_{n-1}^k.
\end{equation}
\end{lemma}

\begin{proof} When the ta-diagram has  only one arc, $a_1$,  the ta-partition is $\{\{1\}\}$,  and  there is only one block,  with  top arc $a_1$.  So,
$$   B_1^1=1,  \quad  A_1=1.$$

We define a  procedure that generates ${{k+3}\choose{3}}$ standard  ta-diagrams  with  $n$  arcs, starting by a standard  ta-diagram  in $\B^k_{n-1}$.     We prove  that every standard ta-diagram in $\AA_n$ is uniquely  generated    by  this  procedure from  one  standard diagram in $\B^k_n$, for some $k$.

Let $\dd$ be a standard  ta-diagram in $\B^k_{n-1}$ with  $k$ blocks.  We  represent a block  by a  black half disc. We put  at  right of $\dd$ the  right  endpoint of the  $n$-th  arc $a_n$.  The left  end point of  this  arc  can  be  put  in $(k+1)$ places,  numbered  by $i=0,1, \dots ,k$; namely: the first at left of $\dd$,  then in $(k-1)$  places between successive blocks, and  finally at right of $\dd$. See Figure \ref{F10}.

For a  given position $i$ of  the left end point of $a_n$,  we consider now the  possible ties form $a_n$  to  the blocks.   When  we  say  that  a  tie  connects $a_n$  to a  block, we  mean  that  the tie  connects $a_n$  to  the top  arc with maximum index of the bloc.   So,  we  can add  at most one  tie form $a_n$ to one of the blocks   at the exterior of $a_n$ (there are $i$ of such blocks); and, contemporarily, we can add at most one tie form $a_n$ to one of the blocks in the interior of $a_n$ (they are $k-i$).  Indeed,  if we  should add  two  ties from $a_n$ to two blocks $\b_r$ and $\b_s$  both  at the exterior  or at the interior of $a_n$,     the  resulting ta-diagram, having two ties connecting a top arc of $\b_r$ and  a top arc of  $\b_s$ to $a_n$,  should be non standard, see  Definition \ref{standard}.

 %Note that these top arcs should result in the  same partition as  $a_n$, and  we want to avoid that, since  among the $(n-1)$ ta-diagrams we  will apply the procedure  to   a  diagram like $\dd$, in which these  two arcs  are tied together, so  already belonging to the same partition. Instead, if $\b_r$ and $\b_s$ are one at the exterior and one at the interior of $a_n$, there exist  no  diagrams like $\dd$, in which they  are connected each other by a  tie.

So, the total  number of standard  ta-diagrams obtained by adding an arc $a_n$  connected  by  ties   to a  diagram $\dd$ with $n-1$ arcs and $k$ blocks,  counting  the  possibilities of  zero  ties at the exterior and  in the interior of $a_n$, is
$$     \sum_{i=0}^k (i+1)(k-i+1)={{k+3}\choose{3}}. $$
Now,  consider any  standard ta-diagram $\dd$ in $\AA_n$.  Since $\dd$ is  standard, the  arc $a_n$    is connected  at most  by one tie  to an arc  at its interior and at most  by one tie  to  an arc in its  exterior. Let's erase these ties and the arc $a_n$. By  erasing $a_n$ and these  ties,  we get a diagram   in $\AA_{n-1}$, for which the sc-partition is obtained from that of $\dd$ by erasing  $n$  from  the part of the sc-partition  of $\dd$ containing $n$.  We  thus   get a  unique standard   ta-diagram in $\B_{n-1}^j $  for  some $j$,   from which  $\dd$  is  uniquely obtained by  the  procedure.

\begin{figure}  [H]
 \centering
 \includegraphics{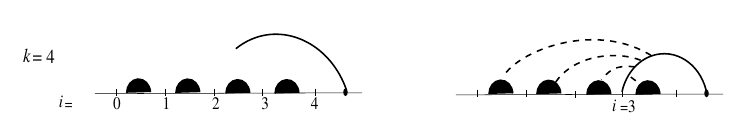}
 \caption{Left: Possible  places for the left endpoints of the arc   $a_n$. Right: Dashed lines  represent  the possibilities  for  placing one tie. }\label{F10}
 \end{figure}
\end{proof}

\begin{remark} \rm The  number  of standard ta-diagrams with $n$ arcs generated by a  standard ta-diagram  $\dd$ with $n-1$  arcs  is  independent from $n$ and  depends only on the  number of  blocks of $\dd$.
\end{remark}

 We  obtain now a stronger  result. In what follows we shall refer to the procedure above simply by the {\it procedure}.

\begin{lemma}\label{L10} \rm
  For  $j=1,\dots,n$,
\begin{equation} \label{recB}  B_n^j= \sum_{k=j-1}^{n-1} B_{n-1}^k {k-j+3 \choose 2}\end{equation}
\end{lemma}

\begin{proof} We  denote  $N_j(\dd^k)$  the number of standard ta-diagrams with $n$ arcs and $j$ blocks generated by the procedure   from one  standard ta-diagram  $\dd^k$ with $n-1$ arcs and $k$  blocs.    Observe that $j$  satisfy  $1\le j\le k+1$.  If $j=k+1$,  there is only one  diagram,  consisting of $\dd^k$ with at right  the  arc $a_n$, without ties,  so  $N_{k+1}(\dd^k)=1$. If $j\le k$, in the new standard diagram with $j$ blocks, say $\dd^j$, there are the first $j-1$ blocks of $\dd^k$, unaltered,  while  the  other  $k-j+1$ blocks of $\dd_k$ will  form  with $a_n$ a sole block. If $k-j=r$, there are $r+2$ possibilities,  since a block $\b$ of
$\dd^k$ belongs to a sole block  of $\dd^j$ containing  $a_{n}$   either if it is  inside  the  new arc $a_n$, or  if it is connected by a tie to $a_n$, or if it is {\it below} the tie from $a_n$  to  the $j$-th  block of $\dd^k$.     See  Figure \ref{F11}.

 \begin{figure}  [H]
  \centering
  \includegraphics{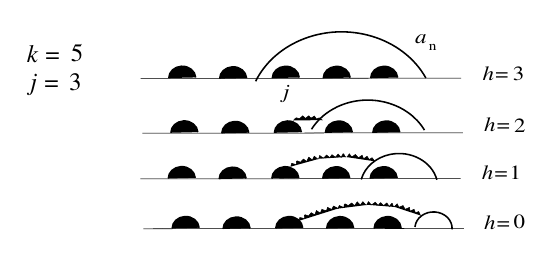}
 \caption{Diagrams with 3 blocs  generated by a  diagram  with 5 blocks. }\label{F11}
  \end{figure}

We must now take into  account  the multiplicities of the above $r+2$ diagrams, by counting the possibilities for adding one allowed tie from $a_n$  to  a  block  inside $a_n$. In fact,  each  diagram of the above list has multiplicity  $1+ h$  where $h$ is the number of blocks inside $a_n$, $h=0,\dots, r+1$. So,
$$  N_j(\dd^k) =  \sum_{h=0}^{k-j+1} (1+h)= {k-j+3 \choose 2}. $$
Therefore Lemma \ref{L10} follows.
 \end{proof}

Comparing   now     Equations (\ref{iniz}) with the  initial conditions (\ref{E1})  and  Equation (\ref{recB}) with  the recurrence  (\ref{E2}) we  conclude that  the  integers $B_n^{n-k}$   fulfill the same  recurrence as   $F(n,k)$,  as states  Equation (\ref{part2}).  Part 2  is proven.

\begin{proof}[Proof  of  Theorem \ref{t3} ]   Because of Equation (\ref{part2}),  observe  that  the three  statements of the  theorem follows  from   the proved lemmas  concerning  ta-diagrams.  More precisely:  statement  (\ref{f3})
follows  from  Equation (\ref{An})  together  with  Equation (\ref{part1}).
    Statement  (\ref{f4}) follows  from
Lemma \ref{L9}, using  (\ref{part1})  and (\ref{formula5}).  Statement  (\ref{f5})  follows from Equation (\ref{part1}).

\end{proof}

\end{document}